\documentclass[12pt]{article}
\usepackage{latexsym}
\usepackage{graphicx}
\newcommand{\al}{\alpha}
\newcommand{\bt}{\beta}
\newcommand{\pal}{\underline{\alpha}}
\newcommand{\pbt}{\underline{\beta}}
\newcommand{\pu}{\underline{u}}
\newcommand{\px}{\underline{x}}
\newcommand{\pv}{\underline{v}}
\newcommand{\pn}{\underline{n}}
\newcommand{\pk}{\underline{k}}
\newcommand{\pl}{\underline{l}}
\newcommand{\pAA}{\underline{{\rm{\huge{\bf A}}}}}
\newcommand{\pBB}{\underline{{\rm{\huge{\bf B}}}}}
\begin{document}
\begin{center}
{\large{\bf Simplified non-Navier-Stokes model of turbulent flow 
and its first numerical realization in D2}}
\vskip .2 in
{\small by Krzysztof Moszynski 
\footnote {University of Warsaw, Institute of Applied Mathematics and Mechanics}}
\footnote {Grant ICM G33-10}
\end{center}
\noindent
{\bf General remarks on the simplified model}
\vskip .1 in
%\noindent
Let us consider the following integro-differential equation
$$\rho_t+\pal\nabla_{\px}\rho-\nu\Delta_{\px}\rho=\kappa\int_{\cal{A}}M(\cdot,\pal,\pbt)d\pbt \leqno(1)$$ 
with properly chosen initial and boundary condition (of Dirichlet type).
This is the equation of the simplified model proposed by Marek Burnat (see [1],[2]) 
with the additional diffusion term 
$-\nu\Delta_{\px}\rho$, where $\nu$ is a (small) positive coefficient.

Similarly as in papers [1] and [2], $t$, $\px$, $\pal$ and $\pbt$ are independent variables, 
\begin{itemize}
\item{$\rho(t,\px,\pal)$ is the {\it $\al$-mass density function},}
\item{$\pal$, $\pbt$ are {\it $\al$-velocities},}
\item{the integral term $\kappa\int_{\cal{A}}M(\cdot,\pal,\pbt)$ is the so called {\it mixer}.}
\end{itemize}
We assume  that
$$t\in [0,T],\ \ T>0,\ \ \px\in\Omega,\ \ \pal,\ \pbt \in\cal{A},$$ 
where $\Omega$ and $\cal{A}$ are rectangles in ${\bf R}^2$. Here we take
$$\Omega=[0,L_1]\times[0,L_2],\ \ {\cal{A}}=[D_1,G_1]\times [D_2,G_2].$$
\vskip .1 in
\noindent
To get usual {\it Euler} quantities it is enough to integrate with respect to the variable $\pal$:
\begin{itemize}
\item{for Euler mass density function 
$$\varrho(t,\px)=\kappa\int_{\cal{A}}\rho(t,\px,\pal)d\pal,$$}
\item{for Euler "impulse density function" (non official terminology!) 
$$p(t,\px)=\kappa\int_{\cal{A}}\pal\rho(t,\px,\pal)d\pal,$$}
\item{for Euler impulse of the mass in $\omega\subset\Omega$
$${\rm imp}(t,\omega)=\int_\omega p(t,\omega)d\px,$$}
\item{for Euler velocity 
$$\pv(t,x)={p(t,\px)\over{\varrho(t,\px)}}=
{{\int_{\cal{A}}\pal\rho(t,\px,\pal)d\pal}\over{\int_{\cal{A}}\rho(t,\px,\pal)d\pal}}.$$}
\end{itemize}
\vskip .5 in
\noindent
We shall define the the function $M(\cdot,\pal,\pbt)$ involved in the {\it mixer}.
Let
$$d=\|\pbt\|\rho(t,\px,\pbt)-\|\pal\|\rho(t,\px,\pal)$$
and 
$$r(d)=
\left\{
\matrix{
{-d\over{1+d}}&{\rm for}&d\geq 0\cr
{-d\over{1-d}}&{\rm for}&d<0\cr
}
\right.,
$$ 
where $\|\cdot\|$ is the Euclidean norm.
Then $M$ is defined by
$$M(\cdot,\pal,\pbt)=
\left\{
\matrix{
r(d)\rho(t,\px,\pal)&{\rm if}&d\geq 0\cr
r(d)\rho(t,\px,\pbt)&{\rm if}&d<0\cr
}
\right.. \leqno(2)
$$
It is easy to verify that
$$\int_{\cal{A}}\int_{\cal{A}}M(\cdot,\pal,\pbt)d\pbt d\pal=0.$$
\vskip .1 in
\noindent
Equation (1) is equivalent to
$$\kappa[\rho_t+{\rm div}_{\px}(\pal\rho)-\nu{\rm div}_{\px}\nabla_{\px}\rho]=
\kappa^2\int_{\cal{A}}M(\cdot,\pal,\pbt)d\pbt.\leqno(3)$$
hence, integrating both sides of (3) with respect to $\pal$ we obtain the following equation for the 
{\it Euler quantities}:
$$\varrho(t,\px)+{\rm div}_{\px}p(t,\px)-\nu{\rm div}_{\px}\nabla_{\px}\varrho(t,\px)=0,\leqno(4)$$
Let us integrate both sides of (4) over $\omega\subset\Omega$.  The Gauss Divergence Theorem
implies the following relation for the mass 
$m(t,\omega)=\int_{\omega}\varrho(t,\px)d\px$
contained in $\omega$:
$$m_t(t,\omega)+\int_{\partial\omega}\pn p(t,\px)dS+
\nu\int_{\partial\omega}\pn\nabla_{\px}\varrho(t,\px)dS=0,\leqno(5)$$
where $\pn$ is the unit vector external normal to the boundary $\partial\omega$.  
\vskip .1 in
\noindent
Observe that equation (5) can be read as the\ \ {\it the Mass Conservation Law} for the model considered.

{\it Any change of the mass in $\omega$ is possible only as the result of fluxes through the boundary
$\partial\omega$}:
\begin{itemize}
\item{of the impulse $\pn p(t,\px)$, and/or}
\item{of the mass    $\pn\nabla_{\px}\varrho(t,\px)$.}
\end{itemize} 
\vskip  1 in
\noindent
{\bf First numerical realization of the simplified model in 2D} 
\vskip .1 in
\noindent
Assume that rectangles $\Omega$ and ${\cal{A}}$ are covered by grids $\Omega_{h}$ and ${\cal{A}}_{ah}$ 
of steps $h_i$ and $ah_i$ $i=1,2$ respectively, where
$$\Omega_h=({x_1 }_{ k_1},{x_2}_{ k_2}),\ {\rm and}\ \ {\cal{A}}_{ah}=({\al_1}_{ l_1},{\al_2}_{ l_2})\leqno (6)$$
where ${x_i}_{k_j}=h_ik_j$ and ${\al_i}_{l_j}=ah_il_j,\ \ \pk=(k_1,k_2),\ \ \pl=(l_1,l_2),$ 
$$h_i={L_i\over M_i},\ \ ah_i={{G_j-D_j}\over{PR_j-MR_j}},$$
$$0\leq k_i\leq M_i,\ \ MR_j\leq l_j\leq PR_j, \ \ i,j=1,2.$$
The grid for time interval  $[0,T], \ \ T>0$, of the time-step $\tau={T\over N}$ is as follows  
$$T_{\tau}=(t_n),\ \ t_n=n\tau,\ \ n=0,1,2,\cdots,N.\leqno(7)$$
\vskip .5 in
For approximation of the function $\rho$ on the grid (6) (7) let us introduce the grid function
$$u_{\pk,\pl}^n=u_{k_1,k_2,l_1,l_2}^n\approx\rho(t_n,\px_{\ k_1,k_2},\pal_{\ l_1,l_2})\leqno(8)$$
where $\px_{\ k_1,k_2}=(h_1k_1,h_2k_2),\ \ \pal_{\ l_1,l_2}=(ah_1l_1,ah_2l_2)$.

\noindent
The grid function defined by (8) has to satisfy the following finite difference equation
$$du_{\pk,\pl}^{n+1}+a_1u_{k_1-1,k_2,\pl}^{n+1}+b_1u_{k_1+1,k_2,\pl}^{n+1}+
a_2u_{k_1,k_2-1,\pl}^{n+1}+b_2u_{k_1,k_2+1,\pl}^{n+1}=\leqno(9)$$
$$=d_1u_{\pk,\pl}^n-a_1u_{k_1-1,k_2,\pl}^n-b_1u_{k_1+1,k_2,\pl}^n-
a_2u_{k_1,k_2-1,\pl}^n-b_2u_{k_1,k_2+1,\pl}^n+$$
$$+{{\tau\kappa}\over2}({\rm{\Huge{\bf F}}}(\pu^{n+1})+{{\rm\Huge{\bf F}}}(\pu^n))+{\bf dir_1}+{\bf dir_2}$$
with

$$a_1=-\lambda_1{{ah_1l_1}\over4}-{{\nu\mu_1}\over2},\ \ b_1=-\lambda_1{{ah_1l_1}\over4}+{{\nu\mu_1}\over2},$$
$$a_2=-\lambda_2{{ah_2l_2}\over4}-{{\nu\mu_2}\over2},\ \ b_2=-\lambda_2{{ah_2l_2}\over4}+{{\nu\mu_2}\over2},$$
$$d=1+\nu(\mu_1+\mu_2),\ \ d_1=1-\nu(\mu_1+\mu_2)$$
$$\lambda_1={\tau\over h_1},\ \ \lambda_2={\tau\over h_2},
\ \ \mu_1={\lambda_1\over h_1},\ \ \mu_2={\lambda_2\over h_2}.$$
\vskip .1 in
\noindent
The argument $\pu^n$ of the function ${\bf F}$ is the long (block) vector 
$${\pu^n}=[\pu^n_{\ \pk,\ \pl_0}|\pu^n_{\ \pk,\ \pl_1}|\cdots|\pu^n_{\ \pk,\ \pl_q}]^T,$$
where $q=(MR_1+PR_1)(MR_2+PR_2)-1$ ($q+1$ is the number of all pairs of indexes $\pl=(l_1,l_2)$).
Function ${\rm{\Huge{\bf F}}}(\cdot)$ is the result of the trapezoidal quadrature in D2 of 
$M(\cdot,\cdot,\cdot)$ over $\cal{A}$ with respect to the variable $\pbt$ (see definition of the mixer). 
Note that ${\rm{\huge{\bf F}}}$ is a term non linear with respect to $\pu$.

\noindent
Dirichlet Boundary conditions are introduced by means of variables ${\bf dir_1}$ and ${\bf dir_2}$:  
$${\bf dir_1}=
\left\{
\matrix{
-a_1(DIRL(n)+DIRL(n+1))&{\rm if}&k_1&=&0\cr
-b_1(DIRR(n)+DIRR(n+1))&{\rm if}&k_1&=&M_1-1\cr
else&0\cr}
\right.,\ \
$$
$$
{\bf dir_2}=
\left\{
\matrix{
-a_2(DIRB(n)+DIRB(n+1))&{\rm if}&k_2&=&0\cr
-b_2(DIRT(n)+DIRT(n+1))&{\rm if}&k_2&=&M_2-1\cr
else&0\cr}
\right..
$$
Here $DIRL(\cdot),\ \ DIRR(\cdot),\ \ DIRB(\cdot),\ \ DIRT(\cdot)$ are Dirichlet conditions
at the left side, right side at the bottom and on the top side of the rectangle $\Omega$ respectively.
Time level is given as argument of $DIR\cdot(\cdot)$.

Equation (9) is simply finite difference approximation of the equation (1) on the grid defined in (6)(7).
Finite difference approximation is obtained as follows:
\begin{itemize}
\item{first partial derivatives with respect to variables $x_1$ and $x_2$:
arithmetic mean of central finite differences at time levels $n$ and $n+1$}
\item{similarly, second partial derivatives with respect to  $x_1$ and $x_2$:
arithmetic mean of second (forward-backward) finite differences at time levels $n$ and $n+1$}
\item{at both time levels $n$ and $n+1$ corresponding Dirichlet boundary conditions are taken into account}
\item{the mixer term is approximated using the trapezoidal 2D quadrature with respect to the variable $\bt$ 
at time levels $n$ and $n+1$; this gives corresponding nonlinear terms $F$. Finally, the arithmetic
mean of levels $n$ and $n+1$ is taken.}
\end{itemize} 
Let $\pAA$ and $\pBB$ be the matrices of dimension 
$$(q+1)M_1M_2\times (q+1)M_1M_2$$
corresponding to the left and right hand side of the linear part of equation (9), respectively.  
We can now write down a compact form of equations (9):
$$\pAA \ {\pu^{n+1}}=
\pBB \ {\pu^n}+{{\tau\kappa}\over2}[{\bf F}(\pu^{n+1})+{\bf F}(\pu^n)]+{\bf DIR}\leqno(10)$$
where ${\bf DIR}$ is the sum of all terms introduced by Dirichlet boundary conditions.
Matrices $\pAA$ and $\pBB$ have a block-diagonal structure, 
$$\pAA=                                                           
\left[                                 
\matrix{
{\bf A}&&&&\cr                            
&{\bf A}&&&\cr
&&{\bf A}&&\cr
\cdot&\cdot&\cdot&\cdot&\cdot\cr
&&&&{\bf A}\cr}
\right],\ \
{\underline{{\rm{\huge{\bf B}}}}}=
\left[
\matrix{
{\bf B}&&&&\cr
&{\bf B}&&&\cr
&&{\bf B}&&\cr
\cdot&\cdot&\cdot&\cdot&\cdot\cr
&&&&{\bf B}\cr}
\right]
$$
each diagonal block corresponds to part of system (9) that depends on a fixed $\pl$ only; hence the block
dimension of the matrices $\pAA$ and $\pBB$ is equal to $q+1\times q+1$.
Matrices ${\bf A}$ and ${\bf B}$ of block dimension $M_2\times M_2$ have also a block structure:
$${\bf A}=
\left[
\matrix{
D&B&&&\cr
A&D&B&&\cr
&A&D&B&\cr
\cdot&\cdot&\cdot&\cdot&\cdot\cr
&&&A&D\cr}
\right],\ \
{\bf B}=
\left[
\matrix{
D_1&-B&&&\cr
-A&D_1&-B&&\cr
&-A&D_1&-B&\cr
\cdot&\cdot&\cdot&\cdot&\cdot\cr
&&&&-B&D_1\cr}
\right]
$$
where $A,\ \ B,\ \ D,\ D_1$ are of dimension $M_1\times M_1$ and
$$D=
\left[
\matrix{
d&b_1&&&\cr
a_1&d&b_1&&\cr
&a_1&d&b_1&\cr
\cdot&\cdot&\cdot&\cdot&\cdot\cr
&&&a_1&d\cr}
\right],\ \
D_1=
\left[
\matrix{
d_1&-b_1&&&\cr
-a_1&d_1&-b_1&&\cr
&-a_1&d_1&-b_1&\cr
\cdot&\cdot&\cdot&\cdot&\cdot\cr
&&&-a_1&d_1\cr}
\right]
$$
$$
A=
\left[
\matrix{
a_2&&&&\cr
&a_2&&&\cr
&&a_2&&\cr
\cdot&\cdot&\cdot&\cdot&\cdot\cr
&&&&a_2\cr}
\right],\ \
B=
\left[
\matrix{
b_2&&&&\cr
&b_2&&&\cr
&&b_2&&\cr
\cdot&\cdot&\cdot&\cdot&\cdot\cr
&&&&b_2\cr}
\right]
$$
\vskip .2 in
\noindent
{\bf Solving system (10)}
\vskip .1 in
\noindent
As is proved in [3], the mixer function $M$ is Lipschitz continuous with 
respect to the variable $\rho$. Hence it follows that function ${\Large\bf F}$
is Lipschitz continuous as well. This fact is very important for method of solving
the nonlinear system (10).   
In order to solve approximately the system (10), one can apply the following simple iteration:
$$\pAA\ {\pu}^{(p+1)}=\pBB\ {\pu}^n+
{{\tau\kappa}\over2}[{\huge{\bf F({\pu}^{(p)})+{\huge{\bf F({\pu}^n)}}]+{\bf DIR}}},\leqno(11)$$
where $p$ is the iteration index. Iteration converges for $\tau$ small enough.
At any time step $n$ we can put ${\pu}^{(0)}={\pu}^n$
and if $\|{\pu}^{(p+1)}-{\pu}^{(p)}\|$ is not too large, then we accept 
${\pu}^{n+1}\approx{\pu}^{(p+1)}$.
At each iteration step it is necessary to solve linear system with the (non-symmetric) matrix 
$\pAA$:
$$\pAA x={\bf f},\leqno(12)$$
where ${\bf f}$ is the sum of all known terms in (11), i.e. terms independent of p+1.
To this end the {\it Richardson iteration} was applied:
$$x_{k+1}=x_k+s r_k,\ \ r_k={\bf f}-{\underline{{\rm{\huge{\bf A}}}}}x_k\leqno(13)$$
with optimally chosen relaxation coefficient $s$. 
Let us estimate the optimal relaxation coefficient $s$ corresponding to maximum norm
$\|[x_1,x_2,\cdots,x_N]^T\|_\infty={\rm max}_{j=1,2,\cdots,N} |x_j|$.
If $x$ is the solution of linear system (12), then for errors $e_k$ and $e_{k+1}$ we have 
$$e_{k+1}=(I-s \pAA)e_k.$$
Hence 
$$\|e_{k+1}\|_\infty\leq\|I-s\pAA\|_\infty \|e_k\|_\infty.$$
Assume 
$$\pAA=
\left[
\matrix{
a_{0,0}&a_{0,1}&a_{0,2}&\cdot&a_{1,N}\cr
\cdot&\cdot&\cdot&\cdot&\cdot\cr
a_{N,0}&a_{N,1}&a_{N,2}&\cdot&a_{N,N}\cr}
\right].
$$
and $a_{i,i}\geq 0$. Then it is easy to verify that
$$\|I-s \pAA\|_\infty=\max_{0\leq i\leq N}(|1-s a_{i,i}|+\sum_{j\neq i}|a_{i,j}|).\leqno(14)$$
If ${\sum_{j\neq i}|a_{i,j}|}\leq a_{i,i}$for all $0\leq i\leq N$, then
${\sum_{j\neq i}|a_{i,j}|}+|1-s a_{i,i}|$ takes  minimal value at $s_{\rm opt}={1\over a_{i,i}}$.
In the case discussed  $a_{i,i}=d>0$ and 
$$\sum_{j\neq i}|a_{i,j}|\leq |a_1|+|a_2|+|b_1|+|b_2|,$$ 
hence
$$\|I-s_{\rm opt}\pAA\|_\infty\leq{{|a_1|+|a_2|+|b_1|+|b_2|}\over d}={\bf NOR}.$$
Looking at the coefficients of equation (9), we infer that the sufficient condition
for convergence of the Richardson iteration, ${\bf NOR<1}$ is satisfied if  proportions
of steps in grids (6)(7) are properly chosen, i.e. if
$\lambda_i={\tau\over h_i}$ and $\mu_i={\lambda_i\over h_i}$ for $i=1,2$ are small enough. 
\vskip .2 in
\noindent
{\bf Computational experiment "Collision"}

\noindent
This experiment was run on the cluster halo2 of the Interdisciplinary Center for Mathematical and 
Computational Modeling of the University of Warsaw. Experiment have to be considered as fully 
"virtual", because coefficients of the model have been taken more or less arbitrarily.
Till now we had no possibility to confront our computational experiments with reality.
In such situation, the results obtained have only a qualitative character.

\noindent
{\bf Description of the experiment.}

\noindent
Four streams of mass enter into the "empty" rectangle $\Omega$ through its four sides.
The entering streams are modeled by the Dirichlet boundary conditions for the function $\rho$,
defined on each of four sides of $\Omega$ by positive functions with graphs of triangular shapes and of 
height increasing in time. Whole experiment contains {\bf 1000 time steps}. After {\bf 120 time steps} the four streams
meet in the center of $\Omega$ and the first period of stagnation begins. But the mass is always entering
into $\Omega$ and some parts of the streams start to go back, rubbing against parts of streams 
moving in opposite directions. After {\bf 520 time steps} the first period of stagnation ends, and first eddies
appear. This period of turbulence ends after {\bf 610 time steps}, and in this moment the second period 
of stagnation begin, which ends after {\bf 870 time steps}. At this point the second period of turbulence starts 
with new eddies appearing.
\vskip .1 in
\noindent
Figures below give certain more interesting stages of this process. The field of unit Euler velocity vectors
on the central part of $\Omega$ can be seen. 
In the preprint version of this paper (see $<www.mimuw.edu.pl/preprints>$) one can find more illustrations
concerning the experiment "Collision".

\noindent
{\bf About the program}

\noindent
The program that realizes the algorithm described above, was build and can be run on the cluster halo2 of the
Interdisciplinary Center for Ma\-the\-ma\-ti\-cal and Computational Modeling of the University of Warsaw.
Ad\-mis\-sib\-le dimension of the problem depends on the number of processors used. 
Coefficients, initial and boundary conditions, as well as the the number of processors can be chosen by the user.
\eject
\noindent
{\bf REFERENCES\\
\begin{enumerate}
\item{M.Burnat, K.Moszy\'{n}ski "On some problems of mathema\-ti\-cal modeling of turbulent flow" J.Tech.Phys., 48,3-4,
171-192, 2007. University of Warsaw, Institute of Applied Mathematics and Mechanics, preprint No 161, 2007.}
\item{M.Burnat, "On some mathematical model of turbulent flow with intensive selfmixing", 
submitted to University of Warsaw, Institute of Applied Mathematics and Mechanics pre\-prints 2011.}
\item{K.Moszynski "On certain numerical application of the time-splitting method" University of Warsaw, Institute of 
Applied Mathematics and Mechanics, preprint No. 201 2011}  
\end{enumerate}}
\eject
\pagestyle{empty}
\noindent
{\bf AFTER 500 TIME STEPS}
\begin{center}
\begin{picture}(4,3)(2.1,5)
\includegraphics[width=600pt]{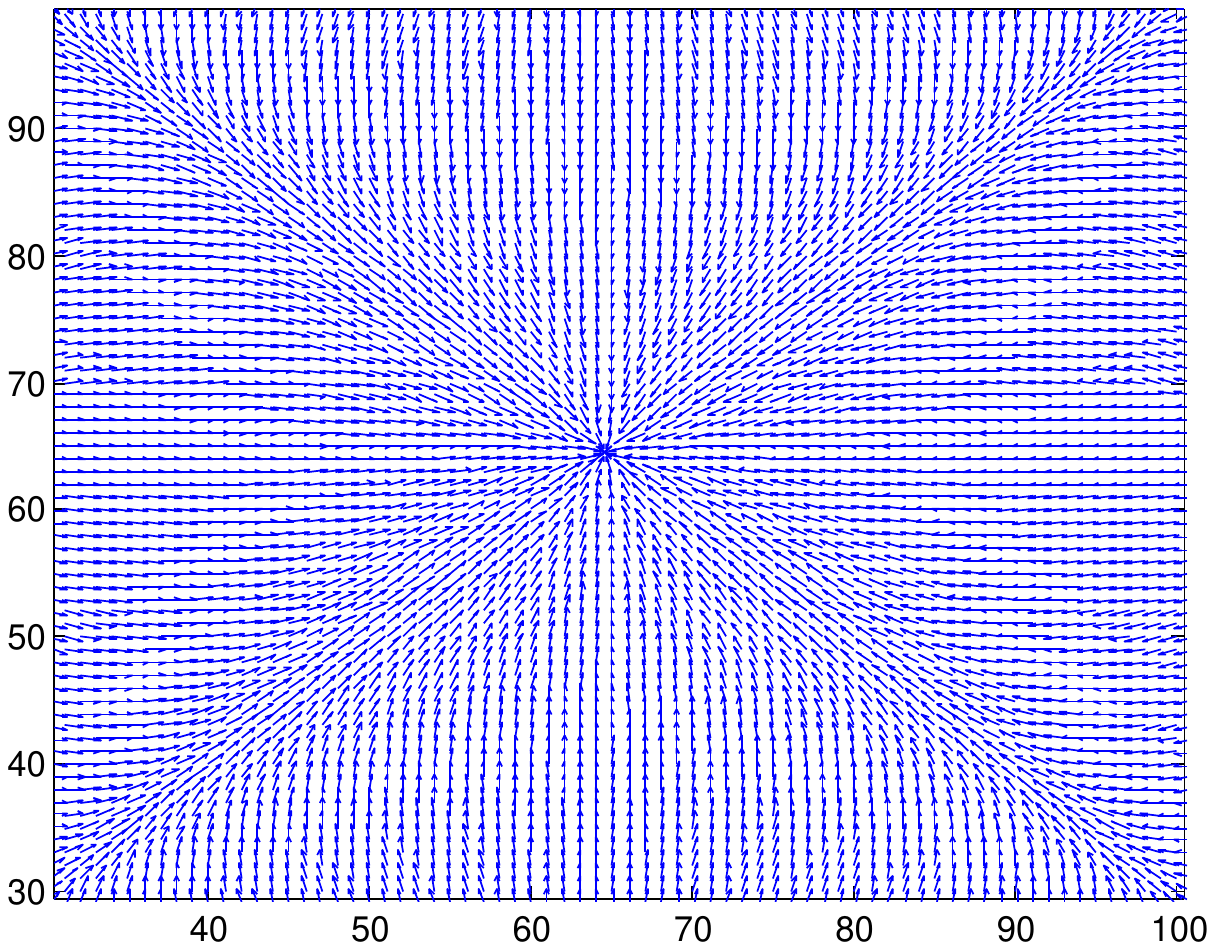}
\end{picture}
\end{center}
\vskip 1 in
\noindent
{\bf AFTER 560 TIME STEPS}
\begin{center}
\begin{picture}(4,3)(2.1,5)
\includegraphics[width=600pt]{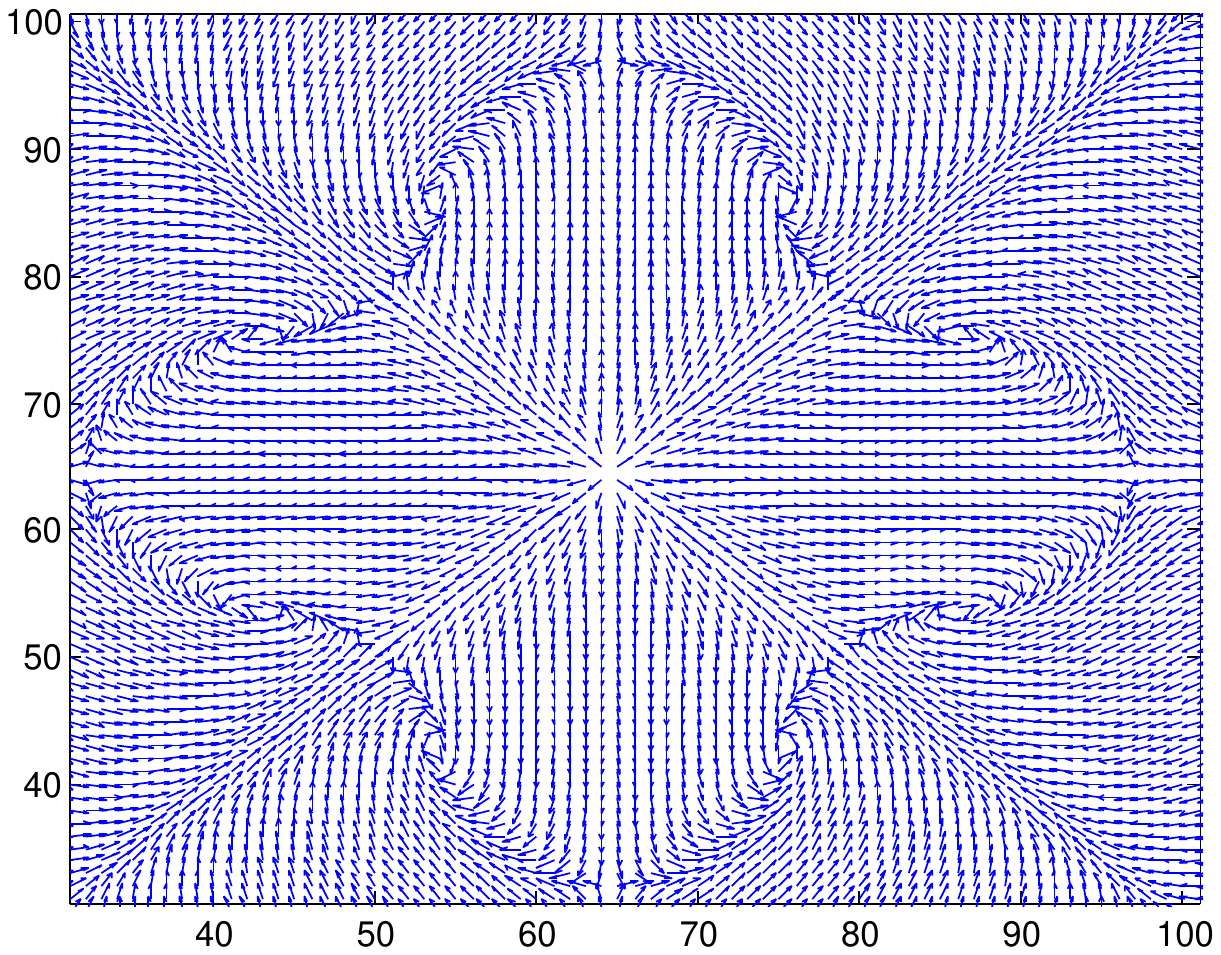}
\end{picture}
\end{center}
\eject
\noindent
{\bf AFTER 850 TIME STEPS}
\begin{center}
\begin{picture}(4,3)(2.1,5)
\includegraphics[width=600pt]{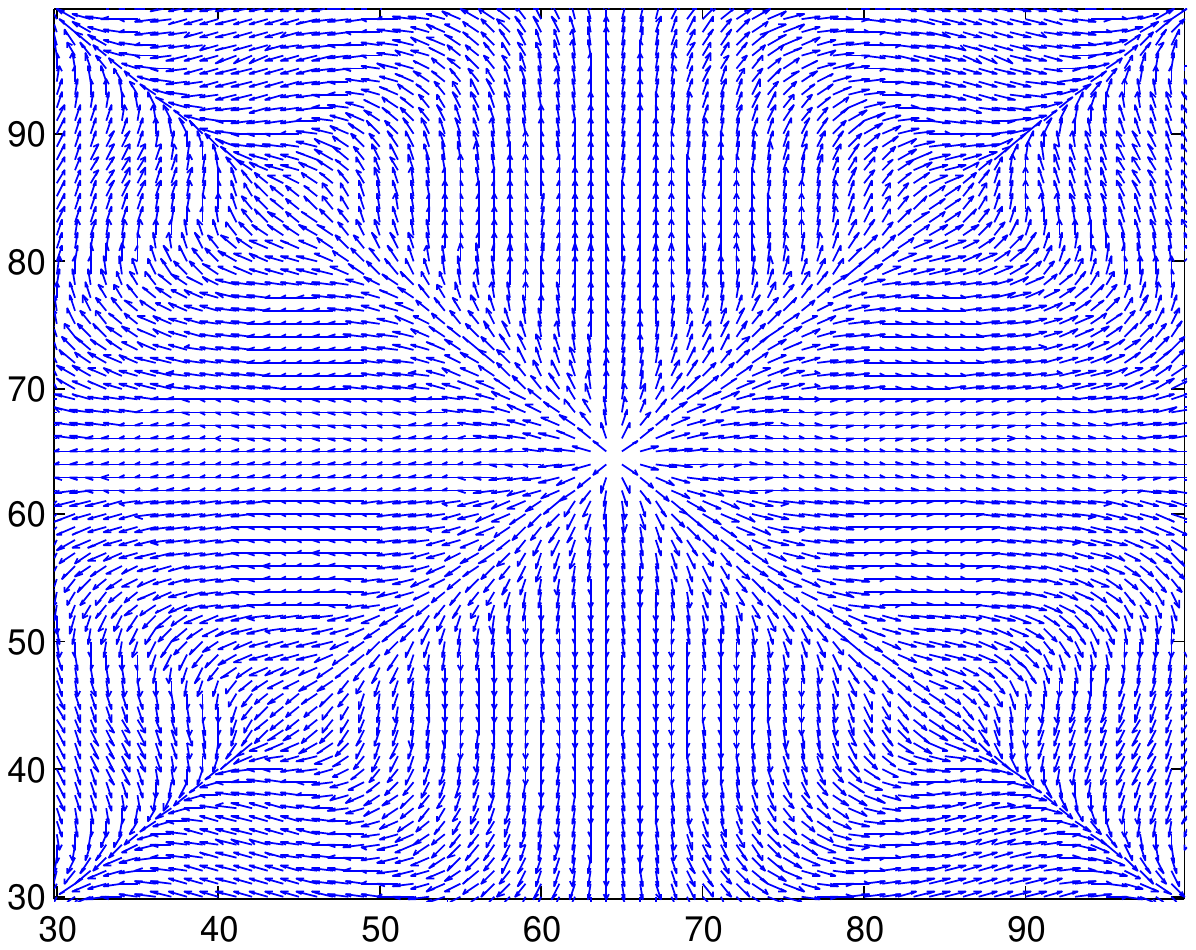}
\end{picture}
\end{center}
\vskip 1 in
\noindent
{\bf AFTER 880 TIME STEPS}
\begin{center}
\begin{picture}(4,3)(2.1,5)
\includegraphics[width=600pt]{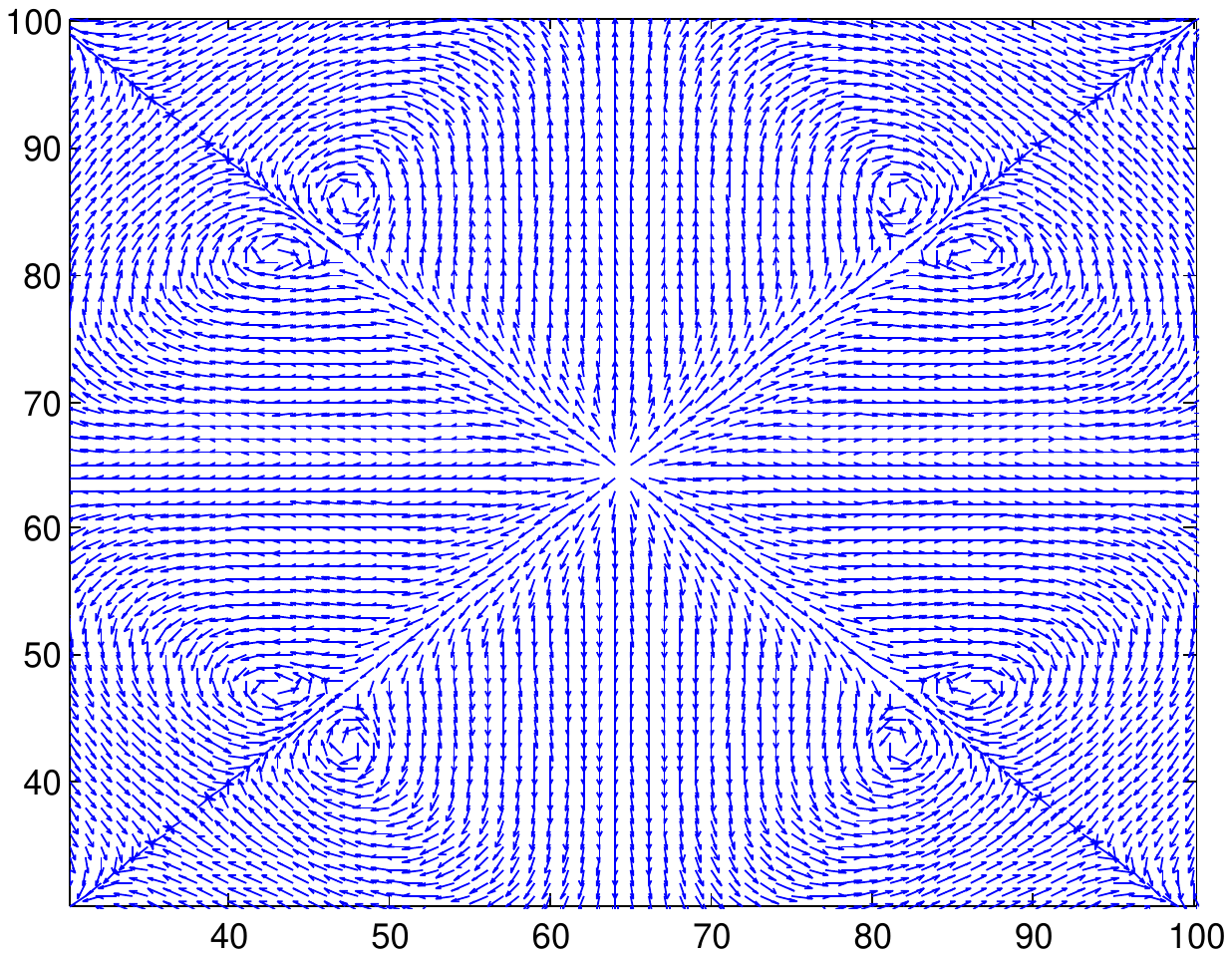}
\end{picture}
\end{center}
\eject
\noindent
{\bf AFTER 1000 TIME STEPS}
\begin{center}
\begin{picture}(4,3)(2.1,5)
\includegraphics[width=600pt]{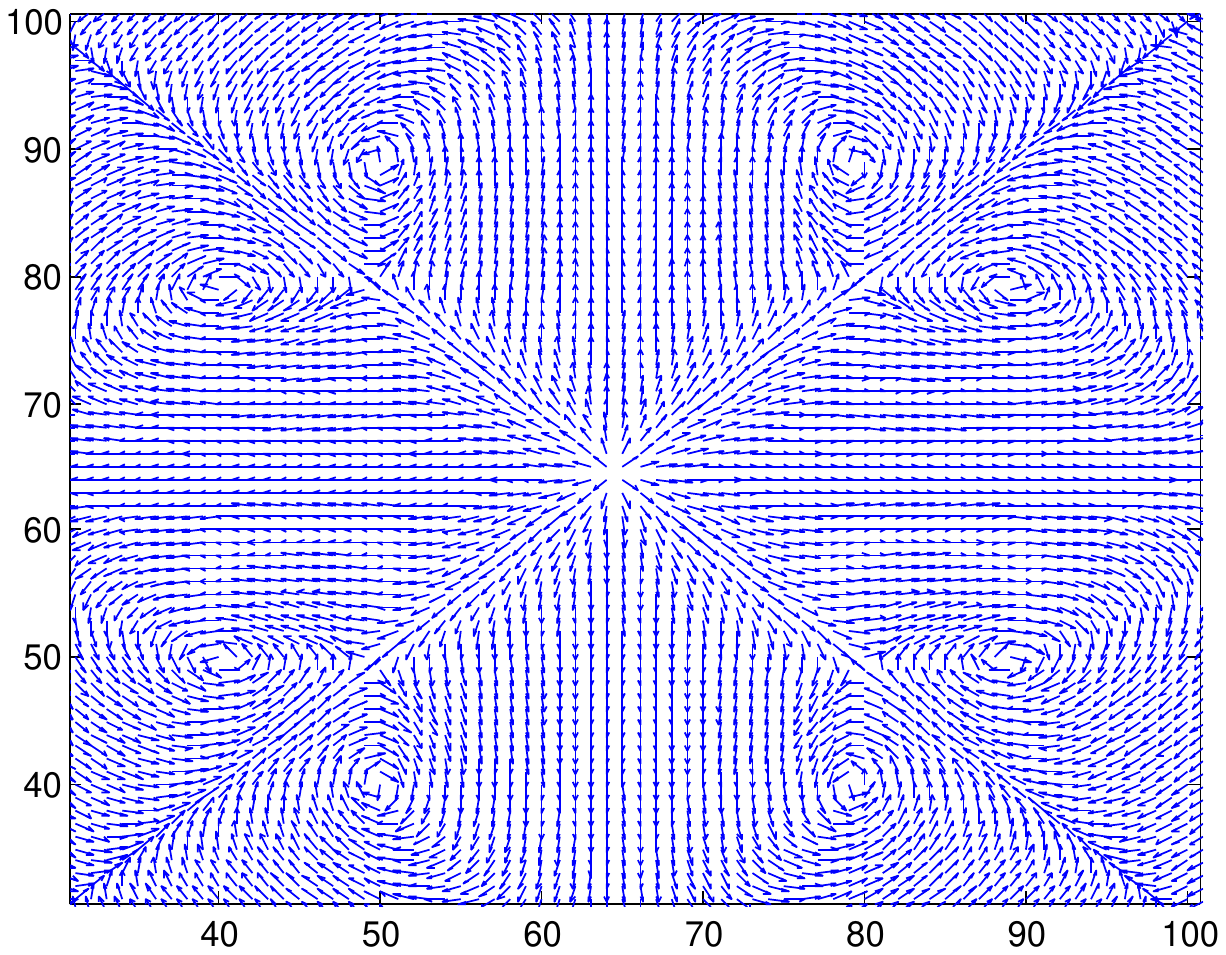}
\end{picture}
\end{center}
\vskip .95 in
\noindent
{\bf AFTER 1000 TIME STEPS (MAGNIFICATION OF A PART)
\begin{center}
\begin{picture}(4,3)(2.1,5)
\includegraphics[width=600pt]{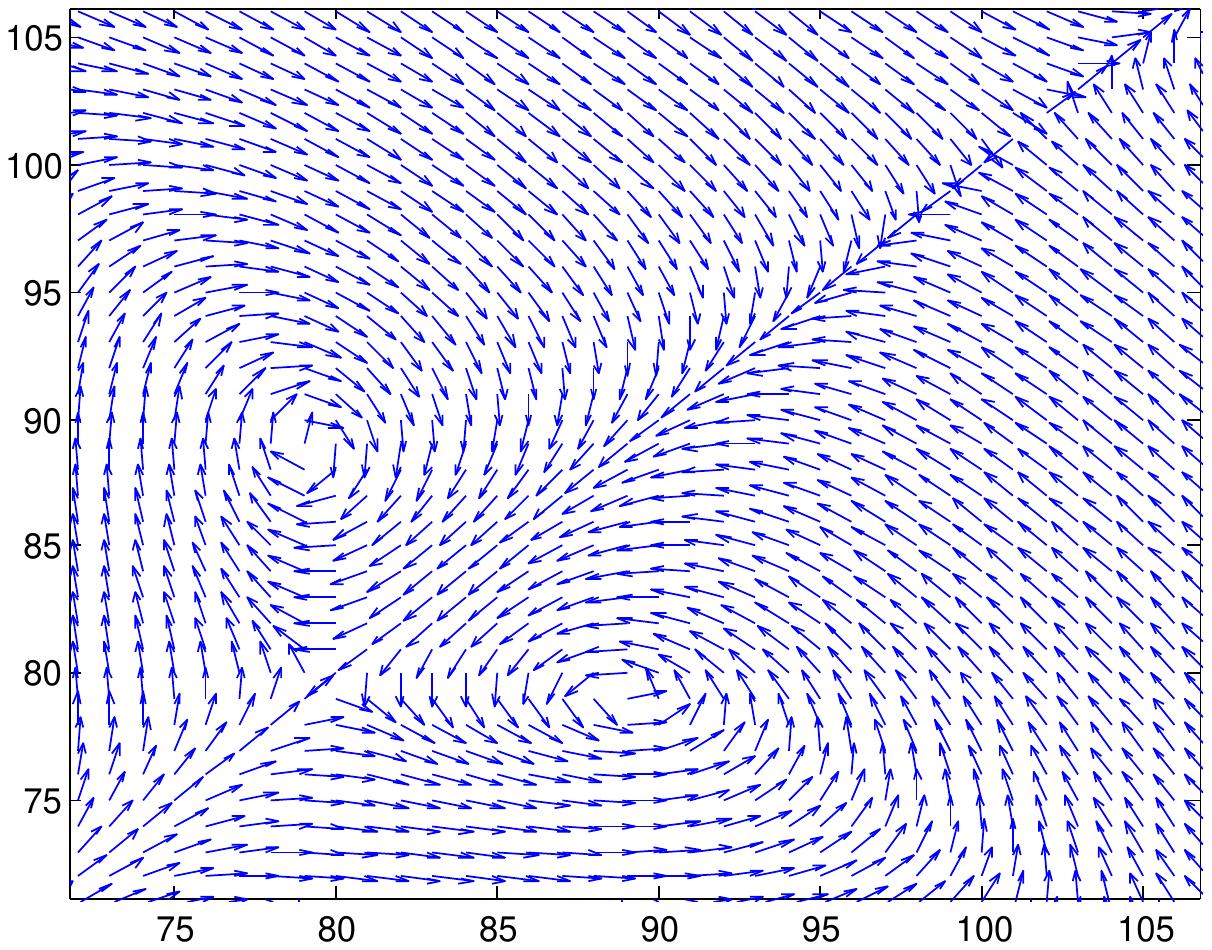}
\end{picture}
\end{center}
\end{document}